# CONFORMAL STRUCTURES AND NECKSIZES OF EMBEDDED CONSTANT MEAN CURVATURE SURFACES

ROB KUSNER

ABSTRACT. Let $\mathcal{M} = \mathcal{M}_{g,k}$ denote the space of properly (Alexandrov) embedded constant mean curvature (CMC) surfaces of genus $g$ with $k$ (labeled) ends, modulo rigid motions, endowed with the real analytic structure described in [15]. Let $\mathcal{P} = \mathcal{P}_{g,k} = \mathcal{R}_{g,k} \times \mathbb{R}_+^k$ be the space of parabolic structures over Riemann surfaces of genus $g$ with $k$ (marked) punctures, the real analytic structure coming from the $3g - 3 + k$ local complex analytic coordinates on the Riemann moduli space $\mathcal{R}_{g,k}$. Then the *parabolic classifying map,* $\Phi : \mathcal{M} \to \mathcal{P}$, which assigns to a CMC surface its induced conformal structure and asymptotic necksizes, is a proper, real analytic map. It follows that $\Phi$ is closed and in particular has closed image. For genus $g = 0$, this can be used to show that every conformal type of multiply punctured Riemann sphere occurs as a CMC surface, and — under a nondegeneracy hypothesis — that $\Phi$ has a well defined (mod 2) degree. This degree vanishes, so generically an *even* number of CMC surfaces realize any given conformal structure and asymptotic necksizes (compare [7, 8] for the case $k = 3$).

## INTRODUCTION

Besides their beauty and variety, perhaps the most important reason minimal surfaces have been so thoroughly investigated is because they can be conformally parametrized via holomorphic functions of a complex variable, a long-familiar tool to many mathematicians. This representation theory was worked out by Enneper and Weierstrass in the middle of the 19th century, leading to important local results and interesting periodic examples. But not until the 1960s did Osserman [20] prepare the way for a global theory by showing that a complete minimal surface with finite total curvature is conformally parametrized by meromorphic data on a finitely punctured, finite genus Riemann surface. The global Enneper-Weierstrass representation enables one to employ algebro-geometric methods from the theory of complex curves to study such minimal surfaces. Over the past two decades it has been used in conjunction with analytic methods by a number of authors — many represented in this volume — to prove deep and striking results about properly embedded minimal surfaces of finite total curvature or finite topology (finite genus and a finite number of ends).

By contrast, the theory of properly embedded surfaces with *nonzero* constant mean curvature (CMC) has developed more recently and without the benefit of these holomorphic methods. Indeed, more than a century elapsed between Delaunay's classification [5] of the CMC rotation surfaces (two ends, genus zero) and Alexandrov's [1] proof that the only compact (zero ends) CMC surface is the round sphere. Only in the past decade have interesting noncompact examples (with three or more ends) of complete CMC surfaces been constructed by Kapouleas [11], using difficult analytic

---

methods. At about the same time, the first steps toward a global theory of these surfaces were taken by Meeks and by Korevaar, Kusner and Solomon: they found topological obstructions to the existence of CMC surfaces with finite topology — none have one end [18], two-ended surfaces are Delaunay unduloids [14] — and developed an asymptotic theory for these surfaces. In fact the fundamental result of [14] is that each (annular) end of a CMC surface has an unduloid asymptote. An immediate consequence is that a finite topology CMC surface can be conformally parametrized by a finitely punctured Riemann surface of finite genus.

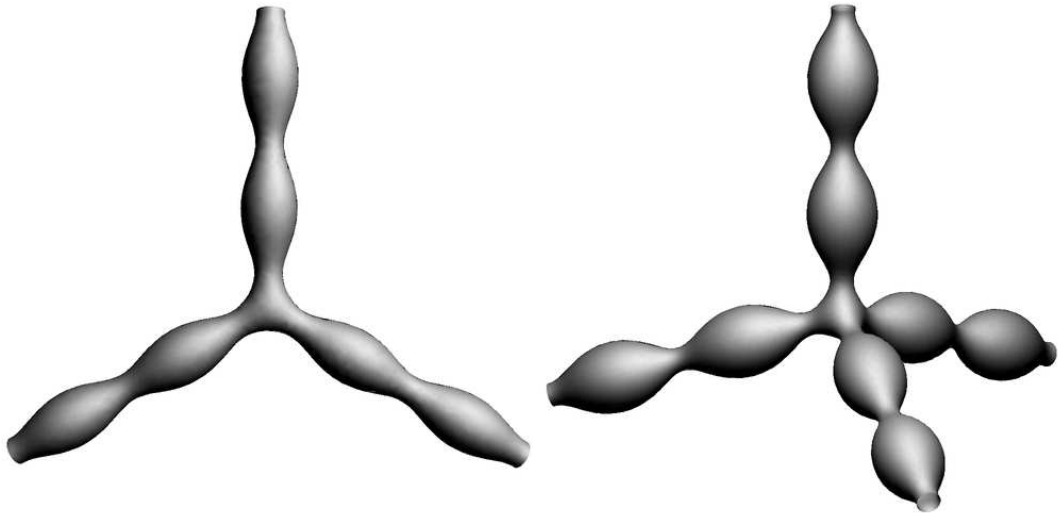

FIGURE 1. Equilateral triunduloid and tetrunduloid (images by Nick Schmitt, GANG)

Of course, this raises the question: Can such a global conformal parametrization be made explicit enough to prove something new, or at least to experiment with CMC surfaces using a computer? Recently at GANG in Amherst (see Figure 1, and www.gang.umass.edu), Schmitt has taken up this question. He adapted a loop-group analogue of the Enneper-Weierstrass representation — dubbed the DPW representation after its discoverers, Dorfmeister, Pedit and Wu [6] — to construct, in terms of explicit meromorphic data on a punctured Riemann surface, conformal immersions of finite topology constant mean curvature surfaces whose ends are asymptotically Delaunay unduloids (or nodoids) [12, 23]. This approach is still under development, and little is known about how to prove the examples it constructs have the symmetries or embeddedness properties required. Nevertheless, it appears to recapture the results of Große-Brauckmann, Kusner and Sullivan [7, 8] about the space of all three-ended CMC surfaces of genus zero, the *triunduloids*.

Our classification of triunduloids made use of transcendental methods, relying heavily on the real analytic variety structure of CMC moduli space developed by Kusner, Mazzeo, and Pollack [15], as well as the compactness results of Korevaar and Kusner [13]. The present paper also uses these transcendental methods, and was motivated by a desire to weave these together with the holomorphic methods arising from the DPW representation to gain a better picture of CMC moduli spaces. Because there is some overlap in these results with those of Mazzeo, Pacard and Pollack, here an effort is made to focus on what is complementary; the reader is encouraged to refer to their nice paper [17], as well as to [8], for additional background and related material.

The author wishes to thank Rafe Mazzeo, Frank Pacard, Dan Pollack, Jesse Ratzkin and Rick Schoen for several stimulating discussions. Thanks are also due to Karsten Große-Brauckmann, Nicos Kapouleas, Hermann Karcher, Nick Korevaar, Franz Pedit, Wayne Rossman, Nick Schmitt, John Sullivan and Mike Wolf for their interest in this work, and especially to David Hoffman and Harold Rosenberg for their invitation to present it at the 2001 Clay/MSRI summer workshop. The research reported here was supported in part by NSF grant DMS-0076085.

## 1. Background on the CMC moduli space

Let $\mathcal{M} = \mathcal{M}_{g,k}$ be the moduli space of all CMC surfaces of genus $g$ with $k$ (labeled) ends, modulo rigid motions. Here, as usual, we scale the mean curvature to be $1$, and include in $\mathcal{M}$ not only the embedded surfaces, but also those which are immersed in the sense of Alexandrov: these are boundaries of immersed 3-dimensional domains (compare [2, 8, 14]). Some of the results mentioned in the introduction can then be expressed in terms of $\mathcal{M}$ as follows:

- $\mathcal{M}_{0,0}$ is a point, represented by any unit sphere, and $\mathcal{M}_{g,0}$ is empty for all $g > 0$ [1];
- $\mathcal{M}_{g,1}$ is empty for every $g$ [18];
- $\mathcal{M}_{0,2}$ consists of the unduloids (up to rigid motion), but $\mathcal{M}_{g,2}$ is empty for all $g > 0$ [14];
- $\mathcal{M}_{g,k}$ is nonempty for every $g \geq 0$, $k > 2$ [11] — for example, $\mathcal{M}_{0,3}$ is a 3-ball [7, 8];
- $\mathcal{M}_{g,k}$ is a real analytic variety of (virtual) dimension $3k - 6$, for $k > 2$ [15].

To understand this real analytic structure on $\mathcal{M}$, it will be useful to introduce the *pre-moduli space* $\widetilde{\mathcal{M}}$, that is, the space of CMC surfaces *before* dividing by rigid motions. We also let $\mathcal{U} = \widetilde{\mathcal{M}}_{0,2}$ denote the pre-moduli space of all unduloids. The discussion below will describe the real analytic variety structure on $\widetilde{\mathcal{M}}$ in terms of an explicit real analytic manifold structure on $\mathcal{U}$, in order to demonstrate the following (almost tautological) result:

**Proposition 1.1.** *The map $\widetilde{\mathcal{M}}_{g,k} \to \mathcal{U}^k$ which assigns a CMC surface its $k$ unduloid asymptotes is real analytic, and is equivariant with respect the actions of rigid motions. Moreover, there is a real analytic necksize function $s : \mathcal{U} \to \mathbb{R}_+$, taking values in $(0, 1]$. It is invariant under rigid motions, and so gives rise to a real analytic asymptotic necksize map $(s_1, \ldots, s_k) : \mathcal{M}_{g,k} \to \mathbb{R}_+^k$, taking values in $(0, 1]^k$.*

*Proof.* The equivariance is clear when we take the diagonal action by rigid motions on the $k$-fold product $\mathcal{U}^k$: as a CMC surface undergoes a rigid motion, each unduloid asymptote undergoes the same motion. Furthermore, once $s$ is defined, it should be clear that $(s_1, \ldots, s_k)$ is defined by composition with the asymptotes map and passage to the quotient by rigid motions.

Recall that an unduloid [14] is a periodic CMC surface of revolution determined by its axis, its neckphase — where along the axis it is nearest the axis — and its necksize. (Of course, modulo rigid motions, an unduloid is specified by its necksize alone.) So the space $\mathcal{U}$ is naturally a fiber bundle over the space $T\mathbb{S}^2$ of (oriented) lines in $\mathbb{R}^3$. The base space $T\mathbb{S}^2$ has a natural analytic structure as a homogeneous space for the group of rigid motions. The total space $\mathcal{U}$ will then be expressed as a 6-dimensional real analytic manifold provided an analytic manifold structure can be exhibited on the typical fiber $\mathcal{D}$, consisting of all unduloids with a fixed, oriented axis (the two ends of the unduloid are labeled according to the axis orientation).

Observe that $\mathcal{D}$ is topologically an open disk, with necksize and neckphase as polar coordinates. The cylinder, whose the necksize is maximal and whose neckphase is not well defined, is at the center of $\mathcal{D}$, fixed under rotation (neckphase-shift). The boundary of $\mathcal{D}$ corresponds to the zero necksize limits of unduloids — chains of spheres along the given axis. The natural analytic structure on $\mathcal{D}$ comes from its (locally homogeneous) realization as an open subset of the unit sphere $\mathbb{S}^2$, indeed as $\mathbb{S}^2 - \{p\}$ for a single point $p$.

This spherical picture of $\mathcal{D}$ also played a key role in [8], leading to a natural interpretation of necksize — there referring to the length of the shortest nontrivial loop on an unduloid — as the spherical distance $n$ between $p$ and another point $q \in \mathcal{D} = \mathbb{S}^2 - \{p\}$. Unfortunately, $n$ is not an analytic coordinate on $\mathcal{D}$: at the cylinder (corresponding to $q$ antipodal to $p$) the function $n$ is only Lipshitz. This will be remedied by reinterpreting *necksize* here to mean $s = \frac{1}{2}(1 - \cos n)$.

Although $s$ may seem geometrically less natural than $n$, it turns out to be more natural in other ways. For example, being a height function on $\mathbb{S}^2$ in the direction $p$, it is symplectically dual to neckphase (angle around $p$) for the homogeneous symplectic structure on $\mathbb{S}^2$. For the present paper, however, the important point is that $s$ is analytic with respect to the natural (locally homogeneous) analytic structure on $\mathcal{U}$, and has the invariance properties and range indicated in the statement of the proposition. In fact, we could have equally well used the *weight* or *force* [14] of an unduloid, $f = n(2\pi - n)$, which is also an analytic coordinate on $\mathcal{U}$, even at the cylinder ($n = \pi$).

The tangent space to $\mathcal{U}$ at an unduloid can be identified with the 6-dimensional linear space $\mathcal{V}$ of *geometric* Jacobi fields (see [14, 15]). Thus $\mathcal{V}$ gives analytic local coordinates on $\mathcal{U}$ in a neighborhood of this unduloid. We shall pretend that $\mathcal{V}$ defines global analytic coordinates on $\mathcal{U}$ in what follows, but since checking analyticity is a local condition, this pretense only amounts to a mild abuse of notation.

Using the asymptotics result [14], one can construct [15] a $6k$-dimensional linear space $\mathcal{W}$ of functions on any surface $\Sigma \in \widetilde{\mathcal{M}}_{g,k}$ which grow at each end of $\Sigma$ like those in $\mathcal{V}$. There are various ways to realize such a $\mathcal{W}$ — either as an orthogonal complement (in an appropriate Hilbert space) to the functions which decay on every end of $\Sigma$, or as a quotient space by such functions — but the salient point is that there is a linear isomorphism from $\mathcal{W}$ to $\mathcal{V}^k$, and thus a real analytic map from $\mathcal{W}$ to $\mathcal{U}^k$.

So in order to establish the proposition above, it would suffice to analytically embed $\widetilde{\mathcal{M}}$ into $\mathcal{W}$, at least locally, in the natural way: at a CMC surface $\Sigma$, represent any nearby CMC surface in $\widetilde{\mathcal{M}}$ via its asymptotic behavior, which is encoded by a point of $\mathcal{W}$. One of the main results in [15] is that this is possible when $\Sigma$ is *nondegenerate*, that is, provided $\Sigma$ supports no nontrivial $L^2$-Jacobi fields. In this case all nearby CMC surfaces form an analytic manifold of dimension $3k$ — half that of $\mathcal{W}$ — and $\widetilde{\mathcal{M}}_{g,k}$ is (locally) a Lagrangian submanifold with respect to a natural symplectic structure defined [15] on $\mathcal{W}$ (or on $\mathcal{V}^k$, or even directly on $\mathcal{U}^k$).

For the general case, a neighborhood of $\Sigma$ in $\widetilde{\mathcal{M}}$ may no longer analytically embed in $\mathcal{W}$, but it can still be embedded as a variety in a larger linear space $\mathcal{W} \times \mathcal{K}$ where $\mathcal{K}$ is a finite-dimensional linear space accounting for the $L^2$-Jacobi fields on $\Sigma$ (compare [15]). Clearly the composition of analytic maps — inclusion of this variety in $\mathcal{W} \times \mathcal{K}$, projection to $\mathcal{W}$, the linear isomorphism to $\mathcal{V}^k$, and finally, the analytic local coordinates map into $\mathcal{U}^k$ — is the asymptotes map. □

## 2. A PROPER CLASSIFYING MAP AND ITS CONSEQUENCES

One can investigate the topology of $\mathcal{M}$ by studying maps with special properties from $\mathcal{M}$ to other known spaces. As noted in the introduction, this is the method used to classify the triunduloids: there is a homeomorphism $\Psi$ from $\mathcal{M}_{0,3}$ to an open 3-ball, realized by the space of ordered triples of points on $\mathbb{S}^2$ (up to rotations) [7, 8]. However, this classifying map $\Psi$ for the triunduloids requires the existence of a symmetry which is not necessarily present for CMC surfaces with four or more ends, and so we will introduce a more general classifying map here.

We have already observed that the *necksize map* $(s_1, \ldots, s_k) : \mathcal{M}_{g,k} \to \mathbb{R}_+^k$, is real analytic. Now consider the *forgetful map* $\phi : \mathcal{M}_{g,k} \to \mathcal{R}_{g,k}$ which assigns the CMC surface $\Sigma$ the conformal class $[\Sigma]$ of its induced metric from $\mathbb{R}^3$. In other words, $\phi$ "forgets" everything but the underlying punctured Riemann surface. It is well known that the Riemann moduli space $\mathcal{R}_{g,k}$ can be given a real analytic structure in two (equivalent) ways, either in terms of finite-area hyperbolic metrics ($6g - 6 + 2k$ real analytic Fenchel-Nielsen parameters) or in terms of holomorphic quadratic differentials with prescribed poles at the punctures ($3g - 3 + k$ complex analytic local coordinates coming from the cotangent bundle to $\mathcal{R}_{g,k}$). Mazzeo, Pacard and Pollack [17] proved a nontrivial result in semilinear elliptic partial differential equations to show the analytic dependence of the unique conformal hyperbolic metric on $\Sigma$, and deduce that the forgetful map $\phi$ is real analytic. (Alternatively, one can exploit the Hopf differential for the Gauss map of a constant mean curvature surface to derive this result via holomorphic quadratic differentials.)

While it is tempting to work with the forgetful map $\phi$ alone, there are at least two reasons why another map will be more suitable. First, in order to develop a degree theory, maps are generally required to be proper. Unfortunately $\phi$ itself is not proper: compare $\mathcal{R}_{0,3}$, which is a single point (and thus compact), to its preimage $\mathcal{M}_{0,3}$ under $\phi$, which is an open (and thus noncompact) 3-ball. In general $\phi$ fails to be proper because a sequence of CMC surfaces diverges in $\mathcal{M}$ by having asymptotic necksizes tend to zero, while remaining in a compact family of conformal types. Second, the meromorphic data for the DPW representation can be geometrically interpreted in terms of flat connections on a rank-2 complex vector bundle over a punctured Riemann surface, with parabolic holonomy determined by a positive number at each of the punctures (compare [3, 23]). Both these considerations suggest replacing the Riemann moduli space with the *parabolic moduli space* $\mathcal{P}_{g,k} = \mathcal{R}_{g,k} \times \mathbb{R}_+^k$ which assigns a positive real number to each puncture of $[\Sigma]$, and defining the *parabolic classifying map*

$$\Phi : \mathcal{M}_{g,k} \to \mathcal{P}_{g,k} = \mathcal{R}_{g,k} \times \mathbb{R}_+^k$$

to be the product of the forgetful map $\phi$ with the asymptotic necksize map $(s_1, \ldots, s_k)$. Clearly $\Phi$ is a real analytic map.

The main result of this paper is the following:

**Theorem 2.1.** *The parabolic classifying map* $\Phi : \mathcal{M} \to \mathcal{P}$ *is proper.*

Before giving a proof of this result (see the next section), we first point out several corollaries:

**Corollary 2.2.** *The parabolic classifying map $\Phi : \mathcal{M} \to \mathcal{P}$ is closed. In fact, its image is a closed real analytic subvariety, which represents a (mod 2) homology class in $\mathcal{P}$. The forgetful map $\phi$ also has closed image in $\mathcal{R}$.*

*Proof.* Since $\mathcal{P}$ is a compactly generated Hausdorff space, we can apply a lemma from general topology (see [8], section 5) to show that any proper map to $\mathcal{P}$ is closed; in particular, $\Phi$ has closed image $\Phi(\mathcal{M}) \subset \mathcal{P}$. The homology class in $\mathcal{P}$ is carried by the image of $\Phi$, or more precisely, by the proper (mod 2) cycle which is the push-forward via $\Phi$ of the fundamental (mod 2) cycle of the real analytic variety $\mathcal{M}$ (compare [4, 24]).

The forgetful map $\phi$ is the composition of $\Phi$ with the projection from $\mathcal{P}$ to $\mathcal{R}$. Unfortunately, projection is generally not a closed map, so the proof that the image of $\phi$ is closed does not follow immediately from the fact that the image of $\Phi$ is closed. Nevertheless, the analyticity (compare [17]) of the forgetful map $\phi$ suffices to show it has closed image. □

For the remainder of this section, assume that the genus $g = 0$. The first observation is that every conformal type of punctured Riemann sphere is realized by a CMC surface. (Although this is also proven quite directly in [15], historical reasons make it interesting to have this argument, based on ideas coming from the original construction of CMC surfaces [11].)

**Corollary 2.3.** *The forgetful map $\phi : \mathcal{M}_{0,k} \to \mathcal{R}_{0,k}$ is surjective.*

*Proof.* From the second part of the preceding corollary, it suffices to show that a dense set of conformal types is realized. Using a little linear algebra, it is not hard to see that, up to conformal transformations of $\mathbb{S}^2$, any configuration of $k > 1$ points $\xi_1, \ldots, \xi_k$ on $\mathbb{S}^2$ can be balanced with *positive* weights: that is, regarding the $\xi_j$ as unit vectors in $\mathbb{R}^3$, there exists a positive solution $f_1, \ldots, f_k$ to the linear relation (force balancing)

$$f_1 \xi_1 + \cdots + f_k \xi_k = 0.$$

For sufficiently small forces $f_j$, the Kapouleas construction [11] then gives a CMC surface $\Sigma$ of genus zero with $k$ ends asymptotic to unduloids whose $j$th axis is approximately in the direction $\xi_j$ and whose corresponding necksize is approximately $f_j$. Because $\Sigma$ is obtained by quasiconformally attaching punctured disks (corresponding to the unduloid ends) at the boundaries of arbitrarily small disks about each $\xi_j$, its conformal structure $[\Sigma]$ can be made to lie within any prescribed neighborhood of $[\mathbb{S}^2 - \{\xi_1, \ldots, \xi_k\}]$ in $\mathcal{R}_{0,k}$. □

It is quite easy to see that $\phi$ cannot be surjective for all $g$ and $k$: when $k = 3$ a CMC surface $\Sigma$ must have a plane of reflection symmetry, so the Riemann surface $[\Sigma] = \phi(\Sigma)$ must have a real conformal involution; this necessary condition, for example, allows only rectangular tori in $\mathcal{M}_{1,3}$.

Next, observe that any component of $\mathcal{M}_{g,k}$ containing a nondegenerate CMC surface is a real analytic variety of dimension $3k - 6$ [15], and this coincides with $6g - 6 + 3k$ — the dimension of $\mathcal{P}_{g,k}$ — precisely when the genus $g = 0$. This is perhaps a third reason to introduce $\mathcal{P}$, since it leads to a degree theory:

**Corollary 2.4.** *For any component $\mathcal{N}$ of $\mathcal{M}_{0,k}$ containing a nondegenerate CMC surface, the restriction of $\Phi$ to $\mathcal{N}$ has a well-defined degree (mod 2). Moreover, this degree is zero.*

*Proof.* Because the domain $\mathcal{N}$ and the range $\mathcal{P}$ have the same dimension, the ideas in [4, 8, 24] let us define this degree as the number of pre-images (mod 2) for any regular value of $\Phi$. The connectedness of $\mathcal{P}$ and properness of $\Phi$ imply that the pre-images of any two regular values can be joined by a compact, 1-dimensional (semi-analytic) variety in $\mathcal{N}$, which necessarily has an even number of endpoints, so the numbers of pre-images agree (mod 2). Then the obvious upper bound on necksize by that of a cylinder shows that this (mod 2)-degree must vanish. □

In particular, it follows that genus zero CMC surfaces generically occur in what might — whimsically and suggestively, considering the proximity of the anticipated neck or bulge phenomena (see [8]) to umbilic points on the surface — be called "innie" and "outie" pairs:

**Corollary 2.5.** *For any nondegenerate surface $\Sigma \in \mathcal{M}_{0,k}$ which is a regular point of $\Phi$, there is a corresponding surface $\Sigma' \in \mathcal{M}_{0,k}$ with the same conformal structure and necksizes.*

## 3. PROOF OF PROPERNESS

The proof that the parabolic classifying map $\Phi : \mathcal{M} \to \mathcal{P}$ is proper relies on the *a priori* estimates for CMC surfaces developed in [13], namely:

- Any $\Sigma \in \mathcal{M}_{g,k}$ lies in a uniform tubular neighborhood of a piecewise-linear graph with $k$ rays and at most $k + 3g - 3$ segments;
- A sequence $\Sigma(i)$ of CMC surfaces, all of which lie in a compact family of uniform tubular neighborhoods, diverges in $\mathcal{M}_{g,k}$ provided the length $\ell(i)$ of the shortest nontrivial loop on $\Sigma(i)$ tends to zero.

A compact family of uniform tubular neighborhoods is characterized by requiring the lengths of all segments to be uniformly bounded above, although these lengths may go to zero, as may the angles between the rays or edges. For sequences of CMC surfaces $\Sigma(i)$ lying in such a compact family of neighborhoods, uniform linear area and total curvature estimates [13] still hold, and uniform pointwise curvature estimates also hold under the hypothesis that the $\ell(i)$ are bounded away from zero (using blow-up arguments of the kind we outline below).

The proof also depends on a description of divergent sequences of punctured Riemann surfaces in terms of the conformal moduli of certain nontrivial annuli. We say that an embedded annulus $A$ on a punctured Riemann surface is *essential* (or nonperipheral) provided it is not homotopic to a single point of the surface, nor to a single puncture. Such an annulus then has a (necessarily finite) conformal modulus $[A] = M \in (1, \infty)$ defined by conformally mapping $A$ into $\mathbb{C}$ so that one boundary component goes to the unit circle, and the other goes to the circle of radius $M$. The basic fact which we will use is the following:

- A sequence of punctured Riemann surfaces $[\Sigma(i)]$ diverges in $\mathcal{R}_{g,k}$ provided there exist essential annuli $A(i) \subset \Sigma(i)$ whose moduli $[A(i)] = M(i) \to \infty$ as $i \to \infty$.

Thus to prove that the classifying map $\Phi$ is proper it will suffice to show that a divergent sequence $\Sigma(i)$ in $\mathcal{M}_{g,k}$, whose asymptotic necksizes $s_1(i), \ldots, s_k(i)$ are uniformly bounded away from zero, must contain divergent essential annuli as above, implying $\Phi(\Sigma(i))$ is divergent in $\mathcal{P}_{g,k}$.

So suppose a sequence $\Sigma(i)$ diverges in $\mathcal{M}$, with all asymptotic necksizes $s_j(i) \geq \ell > 0$. Then we have two alternatives:

- $L(i) \to \infty$, where $L(i)$ is the length of the longest segment in the uniform tubular neighborhood containing an essential subannulus of $\Sigma(i)$; or
- $\ell(i) \to 0$, where $\ell(i)$ is the length of some non-trivial closed curve, which is not homotopic to an end of $\Sigma(i)$.

In the first case, we use the Alexandrov symmetrization method [13, 14] to show there is a sequence of essential annuli $A(i) \subset \Sigma(i)$ having length comparable to $L(i)$ on which they are approximately unduloidal; thus they have moduli $[A(i)]$ growing like $\exp(L(i)) \to \infty$.

In the second case, we make a blow-up argument (compare [8, 13, 14]) to get a non-flat, finite total curvature minimal surface bounding an immersed 3-manifold.

If this blow-up minimal surface is embedded, then its top and bottom ends are catenoidal. Suppose $A$ is a punctured-disk representative of such an end. Then as $i \to \infty$, the end $A$ is approximated by rescaled copies of annuli $A(i) \subset \Sigma(i)$ with moduli $[A(i)] \to [A] = \infty$. Since $A$ is catenoidal, it has non-zero force (see [8, 13, 14]), as must the $A(i)$, which implies they are homologically nontrivial. And although the forces (and thus the necksizes) of the $A(i)$ tend to zero, because we are assuming the necksizes of the ends of $\Sigma(i)$ are bounded away from zero, it follows that the $A(i)$ are not homotopic to any end of $\Sigma(i)$, and so must be essential.

If the blow-up minimal surface is only immersed, the argument is more subtle. If each of its ends is simple, then each end must be either catenoidal or planar. An application of the strong half-space theorem [10], using the immersed 3-manifold in place of $\mathbb{R}^3$, shows that not all ends can be planar, and so we get a catenoidal end, and a punctured-disk representative $A$, to use just as in the argument above.

In case of a non-simple end, its winding number is at least two. Consider a representative $A$ and its curve of intersection with a nearly flat cone in $\mathbb{R}^3$, as well as the intersection of the cone with the immersed 3-manifold. Applying the Gauss-Bonnet formula to the latter shows that this curve of intersection cannot bound in the 3-manifold. It follows that the annuli $A(i) \subset \Sigma(i)$ approximating $A$ must have been homotopically nontrivial in the 3-manifold as well, and thus essential. This completes the proof of our main result.

(The reader may find it amusing to apply this "cone trick" to decide whether a particular immersed surface extends to an immersed 3-manifold. Take, for example, the Enneper surface, a minimally immersed $\mathbb{R}^2$ in $\mathbb{R}^3$; it has one end, which we orient to be asymptotically horizontal, of winding number 3. In general, the inclusion of the bounding minimal surface into the 3-manifold induces an epimorphism on fundamental groups, so if Enneper bounds at all, it must bound an immersed 3-ball with one boundary puncture. Consider the intersection of this immersed 3-ball with a nearly horizontal cone whose vertex lies well "above" the Enneper surface. It gives a nearly flat planar domain with one boundary component, that is, a 2-disk. Applying the Gauss-Bonnet formula to this 2-disk implies that its boundary winds once, not thrice, a contradiction. So Enneper is not immersed in the sense of Alexandrov.)

4. SOME OPEN PROBLEMS

To conclude, it may be interesting to contemplate the following (open) problems which are related to the results above:

- What is the image of the parabolic classifying map $\Phi$ in $\mathcal{P}$, or its intersection with a slice of fixed conformal type? In other words, how does the set of allowable necksizes depend on the underlying punctured Riemann surface? For $g = 0$ and $k = 3$ there is only one conformal type, so these questions coincide: the image of $\Phi$ is known to be the simplex determined by the spherical triangle inequalities for the necksizes of a triunduloid [8]. The general case seems closely related to the Biswas inequalities [3, 23] which arise in the study of flat connections on parabolic vector bundles.
- Is the cardinality of $\Phi^{-1}([\Sigma])$ finite, and is $2^{|\chi(\Sigma)|}$ an upper bound? The absolute value of the Euler characteristic $|\chi(\Sigma)| = 2g - 2 + k$ is the number of trousers in a decomposition of $\Sigma$, so the possible choices of "innie" or "outie" configurations at each of the trousers saturates this bound. Again, for $g = 0$ and $k = 3$ we have $|\chi(\Sigma)| = 1$, and this bound is sharp [8]: $\Phi$ is two-to-one, except at the maximal necksize "fold" where it is one-to-one.
- Is $\mathcal{M}$ connected? And, at least on some connected component of $\mathcal{M}$, how far from an isomorphism is the homomorphism on fundamental groups induced by $\Phi$ (or, equivalently, by $\phi$)? In case $g = 0$ it follows from [8] and [22] (or more directly from [17]) that, on what is presumably the only component — the one containing the *coplanar* CMC surfaces [9] — of $\mathcal{M}_{0,k}$, this is an epimorphism for any $k$.
- Does $\mathcal{M}$ carry a complete metric of nonpositive curvature, in some suitable sense, and is it thus an Eilenberg-MacLane $K(\pi, 1)$-space? The identification of $\mathcal{M}_{0,3}$ with the (hyperbolic) 3-ball in [8] is very suggestive, but it is not clear where to find a metric (hyperbolic or otherwise) directly from the CMC geometry, even in this special case. Note that $\mathcal{P}_{g,k}$ does carry such a metric, and is a $K(\pi, 1)$ with $\pi = \pi_1(\mathcal{R}_{g,k})$, which when $g = 0$ is a (colored) braid group.
- Can one construct a properly (Alexandrov) embedded *minimal* surface of finite topology which is degenerate (in the sense of [16, 19]), either explicitly, or by global methods on minimal surface moduli space (for instance, might one exhibit a classifying map whose critical points must be degenerate minimal surfaces, and which is forced to have a critical point via Morse theory)? And if one works $L^2$-orthogonal to the extra Jacobi fields, can one glue unduloids to the ends of the minimal surface to get a degenerate CMC surface in $\mathcal{M}$? There exist minimal surfaces with all ends planar which *are* degenerate, but these are not immersed in the sense of Alexandrov; similar comments apply to the immersed constant mean curvature tori constructed by Wente [25] and classified by Pinkall and Sterling [21].
- Is there a computable way to detect when an immersed surface in $\mathbb{R}^3$ is immersed in the sense of Alexandrov? As we have noted above, the trick of intersecting with (almost) flat surfaces, such as cones, and applying Gauss-Bonnet seems to do the job for "nice" surfaces, but it would be interesting to have a general obstruction theory for this purely differential topology problem.

CENTER FOR GEOMETRY, ANALYSIS, NUMERICS & GRAPHICS (GANG), DEPARTMENT OF MATHEMATICS, UNIVERSITY OF MASSACHUSETTS, AMHERST, MA 01003
*E-mail address*: kusner@math.umass.edu